\documentclass[12pt]{report}
\usepackage{amssymb,amsmath,makeidx,verbatim}
\newcommand{\C}{\mathbb{C}}

\newcommand{\be}{\begin{enumerate}}
	\newcommand{\ee}{\end{enumerate}}
\newcommand{\bq}{\begin{eqnarray*}}
	\newcommand{\eq}{\end{eqnarray*}}

\begin{document}
	\newcommand{\disp}{\displaystyle}
	\thispagestyle{empty}
	\begin{center}
		\textsc{On the operator-valued Fourier transform of the Harish-Chandra Scwartz algebra\\}
		\ \\
		\textsc{Olufemi O. Oyadare}\\
		\ \\
		Department of Mathematics,\\
		Obafemi Awolowo University,\\
		Ile-Ife, $220005,$ NIGERIA.\\
		\text{E-mail: \textit{femi\_oya@yahoo.com}}\\
		
	\end{center}
	\begin{quote}
		{\bf Abstract.} {\it We establish a $K-$type decomposition of the Harish-Chandra Schwartz algebra $\mathcal{C}^{p}(G),$ for any real-rank $1$ reductive group $G$ with a maximal compact subgroup $K$ and $0<p\leq2.$ This decomposition is then used to give an infinite-matrix-realization of the operator-valued Fourier image $$\mathfrak{F}:\mathcal{C}^{p}(G)\rightarrow \mathcal{C}^{p}(\hat{G})$$ of $\mathcal{C}^{p}(G)$ as a Fr$\acute{e}$chet multiplication algebra in which every member of $\mathcal{C}^{p}(\hat{G})$ consists of a countable block-matrices of the form $$((\mathfrak{F}_{B}(\breve{\alpha})_{(\gamma,m)}(\Lambda)\otimes\mathfrak{F}_{H}(\breve{\alpha})_{(\gamma,l)}(Q:\chi:\nu))_{\gamma\in F, (l,m)\in\mathbb{Z}^{2}})_{F\subset \hat{K},|F|<\infty}$$ for every $\alpha\in \mathcal{C}^{p}(G).$ This proves Trombi's conjecture for $G$ of real rank $1$ and the technique leads to a proof of the fundamental theorem of harmonic analysis for any arbitrary real-rank reductive group $G.$}
	\end{quote}
	\ \\
	\ \\
	\ \\
	\ \\
	\ \\
	\ \\
	\ \\
	\ \\
	\ \\
	\ \\
	\S1. {\bf Preliminaries.}\\
	\indent This paper considers the operator-valued Fourier transform image of the $L^{p}-$ version of {\it Harish-Chandra Schwartz} algebra for any real reductive group $G$ in the {\it Harish-Chandra class.} To this end we introduce a method of reduction of the operator-valued problem to the consideration of only the scalar entries of a countably-infinite blocks of a block of matrices. That is, a reduction to the scalar entries of a countably-infinite matrix consisting of countably-infinite blocks of matrices. This perspective was motivated by Varadarajan $[24.]$ in the presence of some Lemmas on the decomposition of the (necessarily) countable $[27.]$ set $\hat{K}$ of equivalence classes of (necessarily finite) irreducible unitary representations of a maximal compact subgroup $K$ of $G.$
	\ \\
	\indent This feat has been successfully employed for $G=SL(2,\mathbb{R})$ by both Barker $[5.]$ for $0<p\leq2$ as a {\it Fr$\acute{e}$chet} space of linear operators and by Varadarajan $[24.]$ for $p=2$ as a {\it Fr$\acute{e}$chet} algebra of infinite matrices. It is a basic elementary exercise that every linear operator may be easily realized as a matrix consisting of finite or infinite entries. The methods of Barker and Varadarajan, which may together be termed {\it the type-$(m,n)$ approach,} considered $K$ as $$K=\{u_{\theta}=\left(\begin{array}{cc} \cos\theta & \sin\theta \\
		-\sin\theta & \cos\theta
	\end{array}\right):\theta\in(0,2\pi]\}=S0(2)$$ and the finite set $\{m,n\}$ in $\hat{K}=\mathbb{Z},$ using them to discuss (in some details) subjects like the type-$(m,n)$ spherical functions on $SL(2,\mathbb{R}),$ asymptotic growth of matrix coefficients of the irreducible unitary representations of type-$(m,n),$ (constructed) the orthogonal Schwartz-type Fr$\breve{e}$chet algebras $^{o}\mathcal{C}^{p}_{mn}(SL(2,\mathbb{R}))$ and $^{\#}\mathcal{C}^{p}_{mn}(SL(2,\mathbb{R})),$ computed their Fourier transforms, which were put together to arrive at the operator-valued Fourier transform of the entire $\mathcal{C}^{p}_{mn}(SL(2,\mathbb{R})),$ for $0<p\leq2.$
	
	\indent In its details, Varadarajan $[24.]$ used an $L^{2}-$orthogonal basis $\{e_{n}\}_{n\in\mathbb{Z}}$ of the representation space of the regular representation $\pi$ such that $$\pi(u_{\theta})e_{n}=e^{in\theta}e_{n}$$ and which splits into $\pi=(\pi_{k},\pi_{\epsilon,i\mu}):$ the discrete series $\pi_{k\in\mathbb{Z}}$ and the irreducible unitary principal series $\pi_{\epsilon,i\mu\in\{0,1\}\times\mathbb{R}}.$ The type-$(m,n)$ subspace $L^{2}_{mn}(G)$ of $L^{2}(G)$ satisfying the requirement $$f(u_{\theta_{1}}xu_{\theta_{2}})=e^{i(m+n)\theta}f(x),$$ $x\in G,$ was considered (See Ehrenpreis and Mautner $[8.],$ $[9.]$ and $[10.]$), a spanning example of which is the combined subspace of the subspaces of the matrix coefficients of $\pi_{k}$ and $\pi_{\epsilon,i\mu}$ with respect to $\{e_{n}\}$ given respectively as $$\phi_{kmn}(x)=<\pi_{k}(x)e_{m},e_{n}>$$ and $$f_{mn}(\mu:x)=<\pi_{\epsilon,i\mu}(x)e_{m},e_{n}>.$$
	
	\indent Analysis of the perturbed differential equations satisfied by both types of matrix coefficients leads to the understanding of their asymptotic growths. In particular, we have that $$e^{t}f_{mn}(\mu:u:a_{t}:v)\sim p^{+}(\mu:m:n)c^{+}(\mu)e^{i\mu t}+p^{-}(\mu:m:n)c^{-}(\mu)e^{-i\mu t},$$ where $\sim$ means the difference of the two sides is $O(e^{-2t});$ $p^{\pm}$ are polynomials in $\mu,m,n;$ $u,v\in \mathfrak{U}(sL(2,\mathbb{R}))$ (the universal enveloping algebra of the Lie algebra $\mathfrak{g}=sl(2,\mathbb{R})$ of $G=SL(2,\mathbb{R})$) and $c^{\pm}$ are the $SL(2,\mathbb{R})$ $c-$functions. The knowledge of this control of the growth of both $\phi_{kmn}$ and $f_{mn}(\mu:\cdot)$ suggests the consideration of the Harish-Chandra operator-valued Fourier transform on the Schwartz-type subspace $\mathcal{C}^{2}_{mn}(SL(2,\mathbb{R}))$ of $L^{2}_{mn}(SL(2,\mathbb{R}))$ as the map $$f\mapsto\mathfrak{F}(f)=:(\pi_{k}(f),\pi_{\epsilon,i\mu}(f))$$ denoted as $$((\mathfrak{F}_{B}f)(k),(\mathfrak{F}_{L}f)(\mu))$$ and given as $$(\mathfrak{F}_{B}f)(k):=\int_{G}f\bar{\phi}_{kmn}dG$$ and $$(\mathfrak{F}_{L}f)(\mu):=\int_{G}f\bar{f}_{mn}(\mu:\cdot)dG.$$ All these form the scalar-reduction of the theory which involves proving the continuous injection of $f\mapsto\mathfrak{F}_{B}f$ and $f\mapsto\mathfrak{F}_{L}f$ on the direct sum $$^{o}\mathcal{C}^{2}_{mn}(G)\oplus^{*}\mathcal{C}^{2}_{mn}(G)$$ (as contained in $[24.],$ Theorem $20$ and first part of Theorem $57$) and of the exact wave-packet (as contained in $[24.],$ Theorem $33$).
	
	\indent At the end of the long journey on the explicit asymptotic growth estimates, the direct sum $$\mathfrak{F}_{B}(\mathcal{C}^{2}(SL(2,\mathbb{R})))
	\oplus\mathfrak{F}_{L}(\mathcal{C}^{2}(SL(2,\mathbb{R})))$$ was found to be the (infinite) matrix realization of the operator-valued Fourier image of $\mathcal{C}^{2}(SL(2,\mathbb{R}))$ (c.f. Theorem $60$ of $[24.]$). A short {\it coda} on invariant transform theory on $SL(2,\mathbb{R})$ completes the journey, leading to the following characterization.\\
	
	\indent {\bf 1.1 Proposition.} (Varadarajan $[24.]$) If $T$ is a tempered distribution on $G=SL(2,\mathbb{R}),$ then $T$ is invariant iff $T$ vanishes on all commutators, i.e., $T(f*g)=T(g*f),$ for all $f,g,\in\mathcal{C}^{2}(G).$ In this case, there are unique tempered distribution $\hat{T}_{B}$ on $\hat{B}$ and $\hat{T}_{H}$ on $\hat{H},$ with $\hat{T}_{H}$ symmetric, such that $T(f)=\hat{T}_{B}(f_{B})+\hat{T}_{H}(f_{H}).\;\Box$\\
	
	Even though it was never claimed in the work of Varadarajan that the theory being developed in $[24.]$ was based on the {\it \lq finite' $K-$type} results of Trombi $[20.],$ $[21.],$ it is however very clear that the parameters $m$ and $n$ (in $\mathbb{Z}\simeq\hat{K}$) are the $K-$types of members of $\mathcal{C}^{2}_{mn}(SL(2,\mathbb{R}))$ and that the set $\{m,n\}$ was to be considered as a finite subset of the (necessarily countable) set $\hat{K}.$ Hence $F:=\{m,n\}$ is a (representative) finite $K(=SO(2))-$type of $G=SL(2,\mathbb{R})$ and that each $\mathcal{C}^{p}_{mn}(SL(2,\mathbb{R}))$ (in the notation of the $L^{2}(SL(2,\mathbb{R}))$ results of Varadarajan $[24.]$ and the $L^{p}(SL(2,\mathbb{R}))$ results of Barker $[5.]$) may simply be written as $$\mathcal{C}^{p}_{mn}(SL(2,\mathbb{R}))=\mathcal{C}^{p}(SL(2,\mathbb{R}):\{m,n\})=\mathcal{C}^{p}(SL(2,\mathbb{R}):F)_{|_{F=\{m,n\}}},$$ in the notation $\mathcal{C}^{p}(G:F)$ of Trombi $[20.].$
	
	\indent In order to then complete his $L^{p}$ results on $SL(2,\mathbb{R}),$ Barker $[5.]$ appealed to an argument of Arthur which simply says that:\\
	\indent $(i)$ the algebraic (direct) sum of all $\mathcal{C}^{p}_{mn}(SL(2,\mathbb{R}))$ is exactly the subspace $\mathcal{C}^{p}(SL(2,\mathbb{R}),SO(2))$ of $SO(2)-$finite members of $\mathcal{C}^{p}(SL(2,\mathbb{R}))$ and\\
	\indent $(ii)$ the space $\mathcal{C}^{p}(SL(2,\mathbb{R}),SO(2))$ is dense in $\mathcal{C}^{p}(SL(2,\mathbb{R}))$ (See also Warner $[27.],$ Proposition $4.4.3.5$ on p. $256$).\\
	With these facts at hand (especially the denseness result in $(ii)$ above, which is valid for $\mathcal{C}^{p}(G),$ for any real reductive group $G$), Barker only needed to prove his isomorphism theorem for $\mathcal{C}^{p}(SL(2,\mathbb{R}),SO(2)),$ as contained in $[5.],$ Theorem $18.1,$ p. $97,$ where $\mathcal{C}^{p}(SL(2,\mathbb{R}),SO(2))$ is the space of $K-$finite members of $\mathcal{C}^{p}(SL(2,\mathbb{R})).$
	
	\indent It then follows therefore that without any argument of Arthur $[3.],$ as used in Barker $[5.]$ to remove Trombi's finite $K-$type restriction (See $[5.],$ p. $1$), we could easily use the scalar results of $\mathcal{C}^{p}_{mn}(SL(2,\mathbb{R}))$ in Barker $[5.]$ now seen as $\mathcal{C}^{p}(SL(2,\mathbb{R}):F)_{|_{F=\{m,n\}}}$ and is the subspace $$\mathcal{C}^{p}(SL(2,\mathbb{R}):F)_{|_{F=\{m,n\}}}=\mathcal{C}^{p}(SL(2,\mathbb{R}):\{m,n\})
	=\mathcal{C}^{p}_{mn}(SL(2,\mathbb{R}))$$ (in the notation of Varadarajan$[24.]$) to get matrix-realization of members of all $\mathfrak{F}(\mathcal{C}^{p}(SL(2,\mathbb{R})))$ in the style of Varadarajan $[24.].$
	
	\indent The question we set out to answer in this paper is: In what context could Trombi's space $\mathcal{C}^{p}(G:F)$ be seen as the space of type-$F$ ($=\{m,n\}$) spherical function in $\mathcal{C}^{p}(G),$ so that their algebraic direct sum (which is already dense in $\mathcal{C}^{p}(G),$ being the subspaces $\mathcal{C}^{p}(G:F)$) could be deemed sufficiently dense for us to state te isomorphism theorem for the operator-valued Fourier transform of $\mathcal{C}^{p}(G)$. An answer to this question is contained in \S 4.
	
	\indent The present paper is organised as follows. Basic results on the subject (See $[20.]$) are summarised below while a short historical journey to the heart of the problem (in much technical details than the present introduction) is contained in \S 2. Our main results are in \S 3. (where the crucial Lemmas are proved), \S 4. (where we undertake the required $K-$type decompositions and establish Trombi's conjecture for all real-rank $1$ groups $G$) and \S 5. (where we establish the fundamental theorem of harmonic analysis for all real reductive groups). A {\it coda} gives another perspective via $K-$finite functions in \S 6.
	
	\indent Let us give a first detail to our outlook. Let $G$ denote any real reductive group in the Harish-Chandra class (Varadarajan, $[25.]$) with Lie algebra $\mathfrak{g}$ and with a maximal compact subgroup $K$ (in which $rk(G)\geq rk(K)$) so that, whenever $rk(G)=rk(K),$ we choose $B$ for a Cartan subgroup of $G$ contained in $K.$ $\theta$ is a Cartan involution of $G$ fixing $K$ element-wise. If we also denote by $\mathfrak{b}$ the Lie algebra of $B$ then there exists a lattice $L_{B}\subset\mathfrak{b}^{*}_{\mathbb{C}}$ such that $L_{B}$ is isomorphic to the dual $\hat{B}.$ Let $\mathfrak{h}$ denote a $\theta$ stable Cartan subalgebra of $\mathfrak{g}$ with maximal vector part so that $$\mathfrak{g}=\mathfrak{t}\oplus\mathfrak{s}$$ is the corresponding Cartan decomposition of $\mathfrak{g}.$ Set the non-empty set $M$ to be cenralizer$_{K}(\mathfrak{h}\cap\mathfrak{s})$ while $M^{'}$ denotes normalizer$_{K}(\mathfrak{h}\cap\mathfrak{s}).$ Denote $\mathfrak{a}=\mathfrak{h}\cap\mathfrak{s}.$ The set $P(\mathfrak{a})$ (respectively, $P(A)$ with $A:=\exp \mathfrak{a}$) shall denote the set of all those parabolic subalgebras (respectively, parabolic subgroups) whose split rank is $\mathfrak{a}$ (respectively, $A$).
	
	\indent The (finite) Weyl group $W(A):=M^{'}/M$ acts on $\mathfrak{a}^{*}_{\mathbb{C}}$ and on $\hat{M}$ as $$(\mathfrak{w}v)(H)=v(Ad\;\mathfrak{w}^{-1}(H)),\;\;\;H\in\mathfrak{a}^{*}_{\mathbb{C}}$$ and $$(\mathfrak{w}\sigma)(m)=\sigma(\mathfrak{w}^{-1}\;m\;\mathfrak{w}),\;\;\;m\in M,$$ respectively, with $\mathfrak{w}\in W(A).$ Now let $Q=MAN_{Q}$ denote the Langlands decomposition of $Q\in P(A).$ For any $\chi\in \hat{M}$ choose $\sigma\in\chi$ and with $\nu\in\mathfrak{a}^{*}_{\mathbb{C}}$ define $\xi_{\nu}$ on $A$ as $\xi_{\nu}(a)=e^{\nu(\log a)}$ and it follows that $\sigma\otimes\xi_{\nu}$ (defined on $MA$ and) given as $$(\sigma\otimes\xi_{\nu})(ma)=\sigma(m)\xi_{\nu}(a)$$ can be extended to all of $Q$ by making it unit on $N_{Q}.$ The principal series of representations of $G$ is given therefore as $$\pi_{Q,\chi,\nu}=Ind^{G}_{Q}(\sigma\otimes\xi_{\nu}).$$
	
	Consider any $Q\in P(A)$ and for the set $\hat{K}$ consisting of the set of equivalence classes of irreducible unitary (necessarily finite-dimensional) representations of $K,$ we choose $F\subset\hat{K}$ whose cardinality $|F|<\infty.$ We shall denote by $C^{\infty}_{c}(G:F)$ the space of compactly-supported smooth functions on $G$ whose $K-$types are in $F.$ That is, $\alpha\in C^{\infty}_{c}(G)$ satisfying the requirement that $$d(\xi)\int_{K}\bar{\chi}_{\xi}(k)\alpha(\tau_{l}(k)\cdot x)dk=\alpha(x)
	=d(\xi)\int_{K}\alpha(x\cdot\tau_{r}(k))\bar{\chi}_{\xi}(k)dk$$ for all $\xi\in F,$ where $\tau_{l}$ and $\tau_{r}$ denote the left and right actions of $K$ on $G$ and $\chi_{\xi}$ is the character of $\xi.$ Equivalently, $\alpha\in C^{\infty}_{c}(G:F)$ if and only if $$d(\xi)^{2}\int_{K}\int_{K}\bar{\chi}_{\xi}(k_{1})\alpha(\tau_{l}(k_{1})\cdot x\cdot\tau_{r}(k_{2}))\bar{\chi}_{\xi}(k_{2}))dk_{1}dk_{2}=\alpha(x).$$ A relation on $\hat{K}$ may be defined as follows: $\xi_{1}\sim\xi_{2}$ whenever there exists $\alpha\in C^{\infty}_{c}(G:F)$ such that $$d(\xi_{1})d(\xi_{2})\int_{K}\int_{K}\bar{\chi}_{\xi_{1}}(k_{1})\alpha(\tau_{l}(k_{1})\cdot x\cdot\tau_{r}(k_{2}))\bar{\chi}_{\xi_{2}}(k_{2}))dk_{1}dk_{2}=\alpha(x).$$ This is an equivalence relation on $\hat{K},$ whose partitioning of $\hat{K}$ would be needed later. We may set $F(\xi_{j^{*}}):=\{\xi_{j}\in\hat{K}:\xi_{j_{1}}\sim\xi_{j^{*}}\}.$ Clearly $$\hat{K}=\bigcup_{\xi_{j^{*}}} F(\xi_{j^{*}})$$ is a disjoint union.
	
	\indent The Fourier transform of any $\alpha\in C^{\infty}_{c}(G:F)$ can be defined relatively to both the discrete and principal series of representations, $\pi_{\Lambda}$ (where $\Lambda\in L^{+}_{B}$) and $\pi_{Q,\chi,\nu}$ (where $Q\in P(A),$ $\chi\in\hat{M}$ and $\nu\in\mathfrak{a}^{*}_{\mathbb{C}}$), respectively. Indeed, if we denote by $$\alpha\mapsto\mathfrak{F}(\alpha)$$ the operator-valued Fourier transform of $\alpha\in C^{\infty}_{c}(G:F),$ then
	
	$$\mathfrak{F}(\alpha) = \begin{cases}
		\;\;\;(0,\mathfrak{F}_{H}(\alpha)),  \;\;\;\;\;\;\;\;\text{if}\;\;\; rk(G)>rk(K), \\
		(\mathfrak{F}_{B}(\alpha),\mathfrak{F}_{H}(\alpha)),  \;\;\;\text{if}\;\;\; rk(G)=rk(K),
	\end{cases}$$
	where $\mathfrak{F}_{B}(\alpha)(\Lambda):=\pi^{F}_{\Lambda}(\breve{\alpha})$ and $\mathfrak{F}_{H}(\alpha)(\Lambda):=\pi^{F}_{Q,\chi,\nu}(\breve{\alpha}).$ Here we define the map $\breve{\alpha}$ as $\breve{\alpha}(x)=\alpha(x^{-1})$ and we recall the general averaging-operator on $G$ as $$\pi(\alpha)=\int_{G}\alpha(x)\pi(x)dx.$$ Simply put, $\mathfrak{F}(\alpha)=(\mathfrak{F}_{B}(\alpha),\mathfrak{F}_{H}(\alpha))$ (for all $\alpha\in C^{\infty}_{c}(G:F)$) where $\mathfrak{F}_{B}(\alpha)\equiv0$ whenever $rk(G)>rk(K).$ A crucial lemma on the asymptotic growth of $\alpha\rightarrow\mathfrak{F}(\alpha)$ may be deduced as follows from that of $\alpha\rightarrow\mathfrak{F}_{H}(\alpha).$\\
	
	\indent {\bf 1.2 Lemma.} (Trombi $[20.],$ p. $102$) Let $s>0.$ Then given any $N\geq0$ there exists a constant $C_{N}$ such that for all $\alpha\in C^{\infty}_{c}(G:F)$ for which $supp(\alpha)\subseteq\{x\in G:\sigma(x)\leq t(s)\}$ we have that $$||\mathfrak{F}_{H}(\alpha)(Q:\chi:\nu)||\leq C_{N}(1+|\nu|)^{-N}e^{s|\Re \nu|},$$ where the inequality holds for all $Q\in P(A),$ $\chi\in\hat{M}(F)$ and $\nu \in \mathfrak{a}^{*}_{\mathbb{C}}.\;\Box$\\
	
	\indent Now set $x=k\exp(X),$ $k\in K,$ $X\in \mathfrak{s}$ and consider the twin elementary {\it spherical functions} $\Xi(x)=\int_{K}e^{-\rho_{Q}(H_{Q}(x))}dk$ and $\sigma(x)=||X||$ (Trombi $[20.],$ p. $100$). It is clear that $\sigma(xy)\leq\sigma(x)+\sigma(y)$ and $\sigma(x^{-1})=\sigma(x)$ for all $x,y\in G$ (Varadarajan $[24.],$ p. $167$). A collection of continuous seminorms may be placed on $C^{\infty}(G)$ given as $$v^{p}_{a,r}(\alpha)=\sup_{x\in G}|\alpha(x:a)|\Xi^{-2/p}(x)(1+\sigma(x))^{r},$$ for $\alpha\in C^{\infty}(G),$ $a\in\mathfrak{C},$ $r\in\mathbb{R},$ $0<p\leq2;$ through which we then define the $L^{p}-$Schwartz space on $G$ (and of Harish-Chandra type) as $$\mathcal{C}^{p}(G)=\{\alpha\in C^{\infty}(G):v^{p}_{a,r}(\alpha)<\infty,\forall\;a\in\mathfrak{C}, r\in\mathbb{R}\}.$$ It may be shown that $\mathcal{C}^{p}(G)\subset L^{p}(G)$ ($[20.],$ p. $100$ and $[20.],$ p. $231$) and that $C^{\infty}_{c}(G)\subset\mathcal{C}^{p}(G),$ where these inclusions are continuous and with dense images ($[24.],$ p. $341$).
	
	\indent An understanding of these spaces depend on the asymptotic estimate of $\Xi.$ The basic result concerning the asymptotic behaviour of $\Xi$ on $G=SL(2,\mathbb{R})$ is found as follow.\\
	
	\indent {\bf 1.3 Proposition.} (Varadarajan $[24.],$ p. $228$) There are constants $c_{1},c_{2}$ such that (as $t\rightarrow\infty$) $$\Xi(a_{t})e^{t}=c_{1}t+c_{2}+O(t^{4}e^{-3t}).$$ In particular, there is a constant $c>0$ such that $$e^{-t}\leq\Xi(a_{t})\leq ce^{-t}(1+t),\;\;t\geq0.\;\Box$$
	
	\indent The above $SL(2,\mathbb{R})-$result is generalized to all $G$ as follows.\\
	
	\indent {\bf 1.4 Theorem.} (Gangolli and Varadarajan $[11.],$ p. $161$) Let $d$ be the number of reflecting hyperplanes in $W(A)$ (equivalently, the number of short positive roots). Then we can find a constant $C>0$ such that $$e^{-\rho(\log a)}\leq\Xi(a)\leq Ce^{-\rho(\log a)}(1+||\log a||)^{d}$$ for all $a\in cl(A^{+}).\;\Box$\\
	
	\indent With these $\mathcal{C}^{p}(G:F)$ can be shown to be a convolution subalgebra of $\mathcal{C}^{p}(G),$ both of which are Fr$\acute{e}$chet algebras with respect to the collection of (continuous) seminorms above. If $0<p<q\leq2,$ then we have that $$C^{\infty}_{c}(G)\subset\mathcal{C}^{p}(G)\subset\mathcal{C}^{q}(G)\subset\mathcal{C}^{2}(G)$$ which also restricts to $$C^{\infty}_{c}(G:F)\subset\mathcal{C}^{p}(G:F)\subset\mathcal{C}^{q}(G:F)
	\subset\mathcal{C}^{2}(G:F).$$ We shall denote the space of cusp forms in $\mathcal{C}^{2}(G)$ by $^{o}\mathcal{C}^{2}(G)$ ($=\emptyset,$ whenever $rk(G)>rk(K)$), set $^{o}\mathcal{C}^{p}(G)=^{o}\mathcal{C}^{2}(G)\cap\mathcal{C}^{p}(G)$ and define (its orthogonal subspace as) $$\mathcal{C}^{p}_{H}(G)=\{g\in\mathcal{C}^{p}(G):\int_{G}\alpha(x)\bar{g}(x)dx=0,\forall \alpha\in^{o}\mathcal{C}^{p}(G)\}.$$ In close conformity with the definition of $\mathfrak{F}(\alpha)$ above we shall have that
	
	$$\mathcal{C}^{p}(G) = \begin{cases}
		\;\;\;\;\;\;\mathcal{C}^{p}_{H}(G),  \;\;\;\;\;\;\;\;\;\;\;\;\text{if}\;\;\; rk(G)>rk(K), \\
		^{o}\mathcal{C}^{p}(G)\oplus\mathcal{C}^{p}_{H}(G),  \;\;\;\text{if}\;\;\; rk(G)=rk(K).
	\end{cases}$$

   \indent Asymptotic growth of all the members of $\mathcal{C}^{p}(G)$ follows from those given in Lemma $1.2$ and Theorem $1.4.$ The major pioneering result of Arthur (as a consequence of these various estimates and asymptotic growth) is therefore the following.\\
   
   \indent {\bf 1.5 Theorem.} (Arthur $ [3.],$ p. $2.3$ and $3.2$) The map $\mathfrak{F}:f\mapsto\mathfrak{F}(f),$ for $f\in C^{\infty}_{c}(G),$ extends uniquely to an isometry from $L^{2}(G)$ (respectively, from $\mathcal{C}^{2}(G)$) onto $L^{2}(\hat{G})$ (respectively, onto $\mathcal{C}^{2}(\hat{G}).$) If $\mathcal{C}^{2}(G)^{'}$ denotes the topological dual space of $\mathcal{C}^{2}(G),$ then the transpose map $\mathcal{C}^{2}(\hat{G})^{'}\rightarrow\mathcal{C}^{2}(G)^{'}$ is a topological isomorphism$.\;\Box$\\
	
	\indent The above is the $p=2$ version of Trombi's theorem. Our aim in this paper is to extend Trombi's results ($[20.]$) on the operator-valued Fourier  transform $\alpha\mapsto\mathfrak{F}(\alpha)$ (proven only for each $\alpha\in\mathcal{C}^{p}(G:F)$) to include all of $\alpha\in\mathcal{C}^{p}(G)$ (See Lemma $1,$ p. $102$ and p. $85$ of $[20.]$) and to characterize its Fourier image.
	
	\indent The Fourier image of $\alpha\in\mathcal{C}^{p}(G)$ or its subalgebra $\alpha\in\mathcal{C}^{p}(G:F),$ for some values of $p$ and/or some $G,$ had been characterized by several authors notably Ehrenpreis and Mautner ($[8.],\;[9.],\;[10.]$), Langlands ($[15.]$), Arthur ($[1.],\;[2.],\;[3.]$), Harish-Chandra ($[13.]$), Trombi-Varadarajan ($[23.]$), Eguchi ($[7.]$), Trombi ($[20],$ $[21.]$), Barker ($[5.]$), Varadarajan ($[24.]$) and lately Oyadare ($[16.],\;[17.],\;[18.]$). See a synopsis of these results in the next section. The closest of these results to the complete characterization of the operator-valued Fourier transform of all $\mathcal{C}^{p}(G),$ $0<p\leq2,$ as a Fr$\acute{e}$chet convolution Schwartz algebra are those of Trombi ($[20.]$) which he established for the Fr$\acute{e}$chet convolution Schwartz algebra $\mathcal{C}^{p}(G:F)$ in which $F\subset\hat{K}$ with $|F|<\infty.$ We shall show in this paper that, through a decomposition of $\hat{K}$ into a countable union of finite sets (in \S 3.) and using a technique of Varadarajan ($[24.]$ successfully employed for $G=SL(2,\mathbb{R})$), the Fr$\acute{e}$chet convolution Schwartz algebra $\mathcal{C}^{p}(G)$ can be decomposed into a direct sum of the subalgebras $\mathcal{C}^{p}(G:F)$ (in \S 4.), therefore using these decompositions to give the needed complete characterization of the operator-valued Fourier image of $\mathcal{C}^{p}(G),$ here denoted as $\mathcal{C}^{p}(\hat{G}).$ It is then easy to study the structural details of $\mathcal{C}^{p}(\hat{G})$ as a Fr$\acute{e}$chet multiplication algebra.\\
	\ \\
	\S2. {\bf Major Results on Fourier Transform on $G.$}
	
	\indent Historically the three papers of Ehrenpreis and Mautner ($[8.],\;[9.],\;[10.]$) marked the beginning of the construction of the Schwartz function space on a semisimple Lie group as well as the consideration of the Fourier transform of such function space. In these papers only the Lie group $G=SL(2,\mathbb{R})$ (in its Lorentz realization) was considered and its harmonic analysis carried out using the point of view of the representation theory and differential equations introduced by Bargmann ($[4.]$). Thereafter Langlands employed these equations to derive Plancherel related theorems for $SL(2,\mathbb{R})$ in an unpublished manusript of $1962.$ Taking a cue from here, Harish-Chandra ($[15.]$) introduced his Schwartz space of a general semisimple Lie group, which formed the basis for Arthur's $1970$ thesis ($[1.],$ under the guidance of Langlands) and his two other unpublished manuscripts where a serious look is accorded the Schwartz space and its Fourier transform.
	
	\indent Arthur started out with the consideration of the Harish-Chandra Schwartz space $\mathcal{C}^{2}(G)$ built (from $C^{\infty}(G)$) on the $L^{2}-$space of $G$ of real-rank $1$ in $[1.]$ and of any arbitrary real-rank in $[2.]$ and $[3.].$ The real-rank $1$ case was attained by Arthur by first seeking to derive a Schwartz-space Plancherel formula for $G$ ($[1.],$ p. $30-35$) at both $x=1$ (the identity of $G$) and at any $x\in G,$ which then suggested the form and decomposition format of the Fourier image $\mathcal{C}^{2}(\hat{G}),$ while much of the work on the arbitrary real-rank $G$ was done  ($[2.]$ and $[3.]$) via the parametrization of some parabolic subgroups. Arthur's major theorem is that: the map $$\mathfrak{F}:f\mapsto \hat{f}$$ gives a topological isomorphism of $\mathcal{C}^{2}(G)$ onto $\mathcal{C}^{2}(\hat{G})$ (followed by its tempered distribution consequences). It thus became clear from Arthur's papers that harmonic analysis of the Schwartz space of any of such groups $G$ must entail the elucidation of elementary (or $\tau-$) spherical functions, (all) irreducible unitary representations and the asymptotic estimates/growth of their matrix coefficients (via the analysis of their perturbed differential equations).
	
	\indent Arthur did not however consider example by example of $G;$ did not impose any further restriction on $\mathcal{C}^{2}(G);$ did not consider $\mathcal{C}^{p}(G),$ for $p\neq2$ and did not view the $L^{2}-$Schwartz space as a Fr$\acute{e}$chet convolution Schwartz algebra. Nevertheless, the singular fact that $\mathcal{C}^{p}(G)\subseteq\mathcal{C}^{2}(G)$ for all $0<p\leq2,$ makes Arthur's thesis and the two other manuscripts crucial to every other result proved for the $L^{p}-$Schwartz space $\mathcal{C}^{p}(G)$ on $G.$ A major property of these spaces is its algebraic nature and, with the collection of seminorms earlier defined above, we could prove that it is a topological algebra in the following manner.\\
	
	\indent {\bf 2.1 Theorem.} (Gangolli and Varadarajan $[11.],$ p. $255$) The space $\mathcal{C}^{2}(G)$ is an algebra under convolution and $\mathcal{C}^{2}(G//K)$ is a commutative sub-algebra. The map $(f,g)\mapsto f*g$ from $\mathcal{C}^{2}(G)\times\mathcal{C}^{2}(G)$ to $\mathcal{C}^{2}(G)$ is continuous$.\;\Box$\\
	
	\indent At about the same time (in the year $1970$) that Arthur submitted his Yale thesis, P. C. Trombi also submitted a thesis on almost the same subject, but in his own case Trombi considered the general $p-$case of $0<p\leq2$ for the spherical Schwartz spaces $\mathcal{C}^{p}(G//K)$ as a Fr$\acute{e}$chet convolution Schwartz algebra for split rank one $G.$ Having realized that the successful theory of the Schwartz space provides the correct framework for all questions of harmonic analysis on $G,$ Trombi took it as a lifetime project to characterize the entire $\mathcal{C}^{p}(G),$ without any restriction on $p$ or on $G$ or on $\mathcal{C}^{p}(G).$ A first look at the Harish-Chandra Fourier transform, along this line of generalization, is its correct definition and its control by the seminorms on $\mathcal{C}^{p}(G)$ given as follows.\\
	
	\indent {\bf 2.2 Theorem.} (Gangolli and Varadarajan $[11.],$ p. $262$) For any measurable function $f$ satisfying the strong inequality, the integral $$(\mathfrak{F}f)(\lambda)=\int_{G}f(x)\varphi_{\lambda}(x)dx$$ converges absolutely and uniformly for all $\lambda\in \mathfrak{F}_{I}$ and is a continuous function on $\mathfrak{F}_{I}.$ If $r\geq0$ is such that $c:=\int_{G}\Xi^{2}(1+\sigma)^{-r}dx<\infty,$ then $$|\mathfrak{F}f(\lambda)|\leq c v_{1,r}(f),$$ for all $\lambda\in \mathfrak{F}_{I}.\;\Box$\\
	
	\indent This is a compatibility result between the Fourier transform map and members of $\mathcal{C}^{p}(G).$
	
	\indent First the restriction to $G$ with split rank one was jointly removed by Trombi and Varadarajan ($[23.]$) for the spherical Schwartz convolution algebra $\mathcal{C}^{p}(G//K)$ leading to the celebrated Trombi-Varadarajan theorem, without any further restriction on the real reductive group $G$ or on $p.$ We refer the reader to consult Gangolli and Varadarajan ($[11.]$) to see the limit of the spectral theory on $G$ used in the proof of the Trombi-Varadarajan theorem, a first version of whose main result is the following.\\
	
	\indent {\bf 2.3 Theorem.} (Trombi-Varadarajan $[23.],$ p. $273$) Let $(dx,d\lambda)$ be an admissible pair of Haar measures on $G$ and $\mathcal{F}_{I}.$ Then the Harish-Chandra Fourier transform $\mathfrak{F}$ is a topological algebra isomorphism of the convolution algebra $\mathcal{C}^{p}(G//K)$ with the multiplication algebra $\mathcal{S}(\mathfrak{F}_{I})^{\mathfrak{w}}.$ Moreover, we have, for any $f\in \mathcal{C}^{p}(G//K),$ we have $$(\mathfrak{F}f)(\lambda)=\int_{G}f(x)\varphi_{\lambda}(x)dx,\;\;\lambda\in \mathfrak{F}_{I}$$ and $$f(x)=|\mathfrak{w}|^{-1}\int_{\mathfrak{F}_{I}}(\mathfrak{F}f)(\lambda)\varphi_{-\lambda}(x)c(\lambda)^{-1}c(\lambda)^{-1}d\lambda,\;\;x\in G.\;\Box$$\\
	
	\indent If we define $\mathcal{Z}({\mathfrak{F}}^{0}) = \mathcal{S}
	(\mathfrak{a}^*)$ and, for each $\epsilon>0,$ let
	$\mathcal{Z}({\mathfrak{F}}^{\epsilon})$ be the space of all $\C$-valued
	functions $\Phi$ such that  $(i.)$ $\Phi$ is defined and holomorphic
	on $int({\mathfrak{F}}^{\epsilon}),$ and $(ii.)$ for each holomorphic
	differential operator $D$ with polynomial coefficients we have $\sup_{int({\mathfrak{F}}^{\epsilon})}|D\Phi| < \infty.$ The space
	$\mathcal{Z}({\mathfrak{F}}^{\epsilon})$ is converted to a Fr$\acute{e}$chet algebra by equipping it with the
	topology generated by the collection, $\| \cdot \|_{\mathcal{Z}({\mathfrak{F}}^{\epsilon})},$ of seminorms given by $\|\Phi\|_{\mathcal{Z}({\mathfrak{F}}^{\epsilon})} := \sup_{int({\mathfrak{F}}^{\epsilon})}|D\Phi|.$  It is known that $D\Phi$ above extends to
	a continuous function on all of ${\mathfrak{F}}^{\epsilon}$
	([$23.$], pp. $278-279$).  An appropriate subalgebra of
	$\mathcal{Z}({\mathfrak{F}}^{\epsilon})$ for our purpose is the closed
	subalgebra $\bar{\mathcal{Z}}({\mathfrak{F}}^{\epsilon})$ consisting of
	$\mathfrak{w}$-invariant elements of $\mathcal{Z}({\mathfrak{F}}^{\epsilon})$,
	$\epsilon \geq 0.$ The celebrated {\it Trombi-Varadarajan theorem} is stated below.\\
	
	\indent {\bf 2.4 Theorem.} (Trombi-Varadarajan $[23.]$)  Let $0 < p \leq 2$ and
		set $\epsilon = (\frac{2}{p})-1$.  Then the spherical Harish-Chandra
		Fourier transform $f \mapsto \widehat{f}$ is a linear
		topological algebra isomorphism of ${\cal C}^p(G//K)$ onto $\bar{\mathcal{Z}}
		({\mathfrak{F}}^{\epsilon}).\;\;\Box$\\
	
	\indent This theorem was effortlessly lifted to include all the {\it symmetric-space} Schwartz convolution algebra $\mathcal{C}^{p}(G/K)$ by Eguchi ($[7.]$) using the analysis of the spectral theory on $G/K.$ The spectral point of view could however not be pushed further to accommodate all of $\mathcal{C}^{p}(G)$ since the expected spectrum is no longer pure imaginary. (Gangolli and Varadarajan $[11.],$ p. $355$). Trombi then introduced the method of first considering the subalgebra $\mathcal{C}^{p}(G:F)$ consisting of members of $\mathcal{C}^{p}(G)$ that have their $K-$types in a finite subset $F$ of $\hat{K},$ for both the real-rank $1$ real reductive groups in ($[20.]$) and for any arbitrary real-rank real reductive group via $K-$finite functions of ${\cal C}^p(G)$ in an unpublished manuscript (See Barker $[5.]$). There appeared however to be no way to adapt all their methods to the consideration of the Fourier transform to the entirety of $\mathcal{C}^{p}(G).$
	
	\indent In order to get fresh perspectives on the way forward other authors waded in on the subject by making explicit computations of all necessary results for specific groups, notably $G=SL(2,\mathbb{R}).$ It is in this wise should Barker's {\it AMS Memoirs} $[5.]$ on $\mathcal{C}^{p}(SL(2,\mathbb{R}))$ be seen, where arguments in Arthur's $\mathcal{C}^{2}(G)-$case were used to remove the {\it \lq\lq finite-set membership of the $K-$types"} in Trombi's greatest and most successful efforts. There again (in the explicit $\mathcal{C}^{2}(SL(2,\mathbb{R}))-$case) the spectrum of the differential equations that are needed to finally lift Trombi's results failed to be pure-imaginary and seemed irreparable. Another author with another fresh idea for the way forward, as explicitly presented in the $SL(2,\mathbb{R})-$case, was Varadarajan $[24.]$ who surreptitiously weaned the arguments from its over-reliance on {\it \lq\lq finite-set membership of the $K-$types."} In both of these successful attempts on the entire $\mathcal{C}^{p}(SL(2,\mathbb{R}))$ (of Barker $[5.]$) and $\mathcal{C}^{2}(SL(2,\mathbb{R}))$ (of Varadarajan $[24.]$) it was clear that the pathway of first considering type-$(m,n)$ spherical functions in $[24.]$ is the most efficient.
	
	\indent Even though Varadarajan $[24.]$ only considered $\mathcal{C}^{2}(SL(2,\mathbb{R}))$ (out of all of $\mathcal{C}^{p}(SL(2,\mathbb{R}))$), his method of conclusion for the characterization of the operator-valued Fourier transform of the Schwartz convolution algebra is unsuspectingly the closest to completing the works of Trombi in extending the properties of the operator-valued Fourier transform from $\mathcal{C}^{p}(G:F)$ to all of $\mathcal{C}^{p}(G).$ Varadarajan started out with the scalar problem of $\mathcal{C}^{2}_{mn}(SL(2,\mathbb{R}))$ of $L^{2}-$functions $f$ on $G$ satisfying $$f(u_{\theta_{1}}xu_{\theta_{2}})=e^{i(m+n)\theta}f(x)$$ (for all $u_{\theta_{1}},u_{\theta_{2}}\in K=SO(2)$ and $x\in SL(2,\mathbb{R})$) known as the type-$(m,n)$ spherical functions. Noting that $$e^{i(m+n)\theta}f(x)=e^{im\theta}f(x)e^{in\theta}
	=\tau_{m}(u_{\theta_{1}})f(x)\tau_{n}(u_{\theta_{1}}),$$ in which we have used the expression $\tau_{n}(u_{\theta})=e^{in\theta},$ for any $K-$type of $G=SL(2,\mathbb{R})$ (Knapp $[14.],$ p. $627$), we have the transformation relation $$f(u_{\theta_{1}}xu_{\theta_{2}})=\tau_{m}(u_{\theta_{1}})f(x)\tau_{n}(u_{\theta_{1}}).$$ We observe here that, in the language of the general context of Trombi, the finite subset $F$ of $\hat{K}$ is (in the special case of Varadarajan's $SL(2,\mathbb{R})$) simply $F=\{m,n\}.$ This is made clearer if we recall that the map $$f\mapsto\int_{K}\int_{K}\bar{\chi}_{\xi}(k_{1})\alpha(\tau_{l}(k_{1})\cdot x\cdot\tau_{r}(k_{2}))\bar{\chi}_{\xi}(k_{2}))dk_{1}dk_{2}$$ is a continuous map of $C^{\infty}_{c}(G)$ onto the space $C^{\infty}_{c:\xi_{l},\xi_{r}}(G),$ of $\tau=(\tau_{l},\tau_{r})-$ spherical functions on $G,$ given exhaustively as $$\{f\in C^{\infty}_{c}(G): f(\tau_{l}(k_{1})\cdot x\cdot\tau_{r}(k_{2}))=\xi_{l}(k_{1})f(x)\xi_{r}(k_{2}),\forall k_{1},k_{2}\in K,x\in G\}$$ (Barker $[5.],$ p. $9$). Hence, recalling from \S 1. above, we have that $$C^{\infty}_{c:\xi_{l},\xi_{r}}(G)=C^{\infty}_{c}(G:\{\xi_{l},\xi_{r}\})=C^{\infty}_{c}(G:F)$$ in the notation of Trombi. Varadarajan $[24.]$ then considered the parameters of the set $F=\{\xi_{l},\xi_{r}\}$ (written as $\{m,n\}$ for $G=SL(2,\mathbb{R})$) as variables, using them (as matrix-parameters) to compute a matrix realization of every member of the operator algebra $\pi(\mathcal{C}^{2}(SL(2,\mathbb{R}))).$ This is the function-space technique we generalize to the consideration of $\tau-$sperical functions on $G.$
	
	\indent Such functions may be derived via an admissible representation $\pi$ in which $$\pi_{|_{_{K}}}=\sum_{w\in \hat{K}}n_{w}w,$$ so that if $F_{1}$ and $F_{2}$ are finite subsets of $\hat{K}$ and if $\tau_{l}:=\sum_{w\in F_{1}}n_{w}w$ and $\tau_{r}:=\sum_{w\in F_{2}}n_{w}w$ in which $E_{l}$ and $E_{r}$ are the respective orthogonal projections, then $$f(x)=E_{l}\pi(x)E_{r}$$ is a $\tau=(\tau_{l},\tau_{r})-$ spherical function on $G$ and an eigenfunction of $\mathcal{Z}(\mathfrak{g}_{\mathbb{C}}))$ (Knapp $[14.],$ p. $234$). More explicitly, $\tau_{l}(k)u:=\pi_{_{F_{1}}}(k)u$ and $\tau_{r}(k)u:=u\pi_{_{F_{2}}}(k)$ where $\pi_{_{F}}(k):=E_{_{F}}\pi(k)E_{_{F}}$ with $$E_{_{F}}=\sum_{w\in\hat{K}}\pi(\xi_{w}),\;\;\mbox{in which}\;\;\xi_{w}(k):=dim(w)\cdot\chi_{w}(k^{-1}),$$ $\chi_{w}$ being the character of $w.$ (Gangolli and Varadarajan $[11.],$ p. $32$).
	
	\indent Our aim in the present paper is to establish that, for any real reductive group $G,$ we can realize the Fourier image of $\mathcal{C}^{p}(G)$ as a Fr$\acute{e}$chet multiplication algebra of infinite block matrices. A summary of the major results of the present section follows, numbered according to chronology.\\
	\ \\
	
	{\bf Summary of Major Results of the Theory.}\\
	\ \\
	$A.$ \underline{Scalar-valued Fourier tranform on $G$}\\
	\indent For $p=2:$\\
	\indent $(1)$ L. Ehrenpreis and F. Mautner ($1955;$ $1957;$ $1959$): Fourier analysis of $\mathcal{C}_{mn}(SL(2,\mathbb{R}))$ and $\mathcal{C}(SL(2,\mathbb{R}))$ as Schwartz spaces.\\
	\indent $(2)$ R. P. Langlands ($1962$): Fourier analysis of $\mathcal{C}(SL(2,\mathbb{R}))$ as a Schwartz space.\\
	\indent For $0<p\leq2:$\\
	\indent $(3)$ See $B$ below.\\
	\indent $(4)$ P. C. Trombi ($1970$): Fourier analysis of $\mathcal{C}^{p}(G//K),$ for only split rank $1$ groups $G,$ as a Fr$\acute{e}$chet convolution Schwartz algebra.\\
	\indent $(5)$ P. C. Trombi and V. S. Varadarajan ($1971$): Fourier analysis of $\mathcal{C}^{p}(G//K)$ as a Fr$\acute{e}$chet convolution Schwartz algebra.\\
	\indent $(6)$ M. Eguchi ($1979$): Fourier analysis of $\mathcal{C}^{p}(G/K)$ as a Fr$\acute{e}$chet convolution Schwartz algebra.\\
	\indent $(7)$ See $B$ below.\\
	\indent $(8)$ P. C. Trombi ($1982$): Invariant Fourier analysis of $\mathcal{C}^{p}(G:F)$ as a topological vector space.\\
	\indent $(9)$ See $B$ below.\\
	\indent $(10)$ See $B$ below.\\
	\indent $(11)$ O. O. Oyadare ($2019;$ $2022;$ $2023$): Fourier analysis of $\mathcal{C}^{p}(G)$ as a Fr$\acute{e}$chet convolution Schwartz algebra, whose Fourier image is proved to be the scalar subalgebra $$\{(\widehat{\xi_{1}})^{-1}\cdot h\cdot (\widehat{\xi_{1}})^{-1}:h\in\bar{\mathcal{Z}}({\mathfrak{F}}^{\epsilon})\}\subsetneq\mathcal{C}^{p}(\hat{G}),$$ proof of its Bochner theorem and analysis of the canonical wave-packets on $G.$\\
	\indent $(12)$ See $B$ below.\\
	
	\ \\
	$B.$ \underline{Vector-valued Fourier tranform on $G$}\\
	\indent For $p=2:$\\
	\indent $(1)$ See $A$ above.\\
	\indent $(2)$ See $A$ above.\\
	\indent $(3)$ J. Arthur ($1970;$ $1972$): Fourier analysis of $\mathcal{C}(G)$ for real-rank $1$ groups $G$ and arbitrary real-rank groups as a Schwartz space.\\
	\indent $(4)$ See $A$ above.\\
	\indent $(5)$ See $A$ above.\\
	\indent $(6)$ See $A$ above.\\
	\indent $(7)$ See below.\\
	\indent $(8)$ See $A$ above.\\
	\indent $(9)$ See below.\\
	\indent $(10)$ V. S. Varadarajan ($1989$): Fourier analysis of $\mathcal{C}(SL(2,\mathbb{R}))$ as a Fr$\acute{e}$chet convolution Schwartz algebra.\\
	\indent For $0<p\leq2:$\\
	\indent $(7)$ P. C. Trombi ($1981$): Fourier analysis of $\mathcal{C}^{p}(G:F)$ for real-rank $1$ groups and later $\mathcal{C}^{p}(G,K)$ for arbitrary real-rank groups $G$ as a Fr$\acute{e}$chet convolution Schwartz algebra.\\
	\indent $(8)$ See above.\\
	\indent $(9)$ W. H. Barker ($1988$): Fourier analysis of $\mathcal{C}^{p}(SL(2,\mathbb{R}))$ and $\mathcal{C}^{o}(SL(2,\mathbb{R}))$ as Fr$\acute{e}$chet convolution Schwartz algebras.\\
	\indent $(10)$ See above.\\
	\indent $(11)$ See $A$ above.\\
	\indent $(12)$ O. O. Oyadare ($2024$): The present paper.\\
	\ \\
	\ \\
	\ \\
	\ \\
	\ \\
	\ \\
	\S3. {\bf Some Elementary Results on $K-$types.}
	
	\indent The major lessons in the $1955$ to $1969$ works of the masters as reviewed in the earlier sections are that, from Ehrenpreis and Mautner's type-$(m,n)$ to Trombi's finite membership of $K-$types, the $\tau-$spherical functions (in which $\tau=(\tau_{l},\tau_{r})$ as before) are the gateways to accessing the portals of function space theory (Gangolli and Varadarajan $[11.],$ p. $vii$), that $\tau-$spherical functions are intricately looped with the concept of $K-$types (Knapp $[14.],$ p. $215$), that we need not pair up the $K-$types (as done in type$-(m,n)$) and that the entire situation is controlled by the $K-$types of these functions all of which reside in some finite subset of $\hat{K}.$ It may therefore be crucial to closely understand the $K-$types of (spherical) functions in order to move forward from Trombi.
	
	\indent There are however a few distinguished members of $\hat{K}$ that have been the subject of research papers and that have been used to understand the different types of representations of $G$ from algebraic point of view. For these we refer to the relevant chapters of Knapp $[14.]$ as well as the references therein for orientation. Our aim in this section (and as would be applied to prove our main results in the next section) is to exploit the detailed structure of $\hat{K}$ in relation to the understanding of $\alpha\in C^{\infty}_{c}(G)$ satisfying $$d(\xi)^{2}\int_{K}\int_{K}\bar{\chi}_{\xi}(k_{1})\alpha(\tau_{l}(k_{1})\cdot x\cdot\tau_{r}(k_{2}))\bar{\chi}_{\xi}(k_{2}))dk_{1}dk_{2}=\alpha(x),$$ for all $k_{1},k_{2}\in K,$ $x\in G$ and $\xi\in F\subset\hat{K}$ with $|F|<\infty.$
	
	\indent Let us recall that $\hat{K}$ denote the set of equivalence classes of (necessarily finite-dimensional) irreducible unitary representations of $K.$ The following first result on $\hat{K}$ is a bit difficult to notice in the literature.\\
	
	\indent {\bf 3.1 Lemma.} $\hat{K}$ is a countable set. $\Box$\\
	
	\indent We refer to Warner $[27.],$ p. $250$ for the proof of Lemma $3.1$ and note that a careful look at the detailed verification of Peter-Weyl theorem suggests a knowledge of Lemma $3.1.$ See also Wallach $[26.],$ p. $36$ and Dixmier $[6.],$ \S 15. The content of this Lemma is that the set $\hat{K}$ is equivalent to any one (or to a finite or a countable union) of countable sets like $\mathbb{N},\cdots,$ $\mathbb{Z},\cdots$ $\mathbb{Q},\cdots$ of the field of number theory.\\
	
	\indent {\bf 3.2 Examples.} $(1)$ For $G=SL(2,\mathbb{R})$ we have $K=SO(2).$ Here $\hat{K}\simeq\mathbb{Z}$ where $n\in\mathbb{Z}$ corresponds to the representation $$\tau_{n}:u_{\theta}\mapsto\tau_{n}(u_{\theta})=e^{in\theta}$$ (Knapp $[14.],$ p. $627$).
	
	\indent $(2)$ For $G=SL(2,\mathbb{C})$ we have $K=SU(2)$ and $\hat{K}=\{n\in\mathbb{Z}:n\geq0\}.$\\
	
	\indent Lemma $3.1$ also allows us to use the elementary properties of a countable set when employing the set $\hat{K}.$ One of such properties is the elementary fact that every countable set is a countable (respectively, finite) union of some finite (respectively, countable) subsets. In the special case of $\hat{K}$ we have the following.\\
	
	\indent {\bf 3.3 Corollary.} Let $K$ and $\hat{K}$ be as above. Then there exists $F_{\alpha}\subseteq\hat{K},$ where $\alpha$ is in some countable indexing set $\Omega$ and $|F|<\infty,$ such that $$\hat{K}=\bigcup_{\alpha\in\Omega}F_{\alpha}.$$
	
	\indent {\bf Proof.} We already have from \S 1 that $$\hat{K}=\bigcup_{\xi_{j^{*}}} F(\xi_{j^{*}}),$$ where $\xi_{j^{*}}$ is in some subset of $\hat{K}.$ Lemma $3.1$ simply assures that this subset of $\hat{K}$ is indeed countable$.\;\Box$\\
	
	\indent The finite set $F_{\alpha}$ may be taken as a singleton or as $F_{\alpha}=\{m_{\alpha},n_{\alpha}\},$ etc. We have already seen the case of $F_{\alpha}=\{m_{\alpha},n_{\alpha}\}$ in the type-$(m_{\alpha},n_{\alpha})$ point of view as it was unconsciously undertaken in Ehrenpreis and Mautner ($[8.]$), Barker ($[5.]$) and Varadarajan ($[24.]$). We may specifically refer to the projections $$\mathcal{C}^{2}(G)\rightarrow\mathcal{C}^{2}_{mn}(G)$$ in the proof of Theorem $57$ of $[24.]$ and of Proposition $5.1$ in $[5.].$ The direct import of these elementary results on $\mathcal{C}^{p}(G)$ is our major interest.\\
	\ \\
	
	\indent {\bf 3.4 Proposition.} Every $K-$type of members of $\mathcal{C}^{p}(G)$ is in a finite subset of $\hat{K}.$
	
	\indent {\bf Proof.} Since the projections $$\mathcal{C}^{p}(G)\rightarrow\mathcal{C}^{p}(G:F)$$ is continuous via the map $$\alpha(x)\mapsto d(\xi)^{2}\int_{K}\int_{K}\bar{\chi}_{\xi}(k_{1})\alpha(\tau_{l}(k_{1})\cdot x\cdot\tau_{r}(k_{2}))\bar{\chi}_{\xi}(k_{2}))dk_{1}dk_{2},$$ the result follows from Corollary $3.3.\;\Box$\\
	\ \\
	\S4. {\bf $K-$type Decomposition of $\mathcal{C}^{p}(G)$ and $\mathcal{C}^{p}(\hat{G})$ for real-rank $1.$}
	
	\indent The motivating result of this paper as proved in the last section is the following.\\
	
	\indent {\bf 4.1 Lemma.} For any $G$ of arbitrary real-rank, $$\mathcal{C}^{p}_{\xi_{l},\xi_{r}}(G)=\mathcal{C}^{p}(G:F)_{|_{F=\{\xi_{l},\xi_{r}\}}}.\;\Box$$
	
	\indent Now let $G$ be or real-rank $1.$ We need the space $\mathcal{C}^{p}(G:F),$ for $F\subset\hat{K}$ with $|F|<\infty.$ It has been shown that $\mathcal{C}^{p}(G:F)$ is a Fr$\acute{e}$chet convolution Schwartz algebra which is explicitly seen as $$\mathcal{C}^{p}(G:F) = \begin{cases}
		\;\;\;\;\;\;\;\;\mathcal{C}^{p}_{H}(G:F),  \;\;\;\;\;\;\;\;\;\;\;\;\;\;\;\;\text{if}\;\;\; rk(G)>rk(K), \\
		\mathcal{C}^{p}_{B}(G:F)\oplus\mathcal{C}^{p}_{H}(G:F),  \;\;\;\;\;\text{if}\;\;\; rk(G)=rk(K).
	\end{cases}$$
	
	We refer to \S1 and \S2 for some unexplained notations and to Trombi $[20.],$ p. $108$ for the meanings of the spaces $\mathcal{C}^{p}_{B}(\hat{G}:F)$ and $\mathcal{C}^{p}_{H}(\hat{G}:F).$
	
	\indent The main result of Trombi $[20.]$ is that the operator-valued Fourier transform $\mathfrak{F}$ (as defined in \S 1) is a linear topological algebra isomorphism of $\mathcal{C}^{p}(G:F)$ (under convolution) onto $\mathcal{C}^{p}(\hat{G}:F),$ where $\mathcal{C}^{p}(\hat{G}:F)$ is the Schwartz space of matrix-functions (under multiplication) in which $$\mathcal{C}^{p}(\hat{G}:F) = \begin{cases}
		\;\;\;\;\;\;\;\;\mathcal{C}^{p}_{H}(\hat{G}:F),  \;\;\;\;\;\;\;\;\;\;\;\;\;\;\;\;\text{if}\;\;\; rk(G)>rk(K), \\
		\mathcal{C}^{p}_{B}(\hat{G}:F)\times\mathcal{C}^{p}_{H}(\hat{G}:F),  \;\;\;\;\;\text{if}\;\;\; rk(G)=rk(K),
	\end{cases}$$ which are linear subspaces of functions $g=(g_{B},g_{H})$ for which $$g_{B}(\Lambda,\Lambda,M)=\sum_{i^{'}\in I^{'}_{P}}\mathcal{F}(\alpha_{i^{'}})(\Lambda,\Lambda,M)g_{H}(i^{'})$$ for all $\alpha_{i^{'}}\in\mathcal{C}^{p}(G:F)$ and $\Lambda\in L^{+}_{B}$ such that $$|\Lambda(H_{\beta})|\leq(\frac{2}{p}-1)k(\beta)$$ (for some positive non-compact root $\beta$ of $(\mathfrak{g},\mathfrak{b})$ with $k(\beta):=\frac{1}{2}\sum_{\alpha\in P}|\alpha(\bar{H}_{\beta})|$ where $P$ is a positive root system for $\triangle(\mathfrak{g},\mathfrak{b})$ for which $\beta(\bar{H}_{\beta})=2$). It is noted here that $g_{B}(\Lambda,\Lambda,M)$ denotes the matrix representation of the operator $g_{B}(\Lambda)$ relative to the basis $$\{\phi_{\gamma,l}(\Lambda):\gamma\in F, 1\leq l\leq n(\Lambda:\gamma)\}$$ while the matrix $g_{H}(i^{'})$ is relative to the basis $$\{\phi_{\gamma,m}:\gamma\in F, 1\leq l\leq n(\chi:\gamma)\}.$$ The induced topology from the tensor product topology on $$\mathcal{C}^{p}_{B}(\hat{G}:F)\otimes\mathcal{C}^{p}_{H}(\hat{G}:F)$$ converts $\mathcal{C}^{p}(\hat{G}:F)$ into a Schwartz Fr$\acute{e}$chet space and, by virtue of the Fourier map, $$\mathfrak{F}:\mathcal{C}^{p}(G:F)\rightarrow\mathcal{C}^{p}(\hat{G}:F),$$ we conclude that $\mathcal{C}^{p}(\hat{G}:F)$ is therefore a Schwartz Fr$\acute{e}$chet multiplication algebra.\\
	
	\indent {\bf 4.2 Proposition.} Let $F\subset\hat{K}$ with $|F|<\infty.$ Then $$\bigoplus_{F\subset\hat{K}}\mathcal{C}^{p}(G:F)$$ is dense in $\mathcal{C}^{p}(G).$
	
	\indent {\bf Proof.} We employ Theorem $4.4.2.1$ (on p. $261$) and Theorem $4.4.2.5$ (on p. $265$) of Warner $[27.]$ and the Remark of Warner $[27.],$ p. $161$ to conclude that for every $\alpha\in\mathcal{C}^{p}(G),$ the sum $$\sum_{F\subset\hat{K},\\|F|<\infty}\bar{\chi}_{_{F}}*\alpha*\bar{\chi}_{_{F}}$$ converges absolutely to $\alpha.$ Hence the proof$.\;\Box$\\
	
	\indent {\bf 4.3 Remarks.}\\
	\indent $(1)$ Proposition $4.2$ above provides a technical way of handling any problem or gap that could be created by the existence of functions $\alpha\in\mathcal{C}^{p}(G)$ with a $K-$type not in a finite subset of $\hat{K},$ if such functions ever exist (See Proposition $3.4$). It simply says that any of such functions (if they exist) could be taken care of by a countable sum of members of $\mathcal{C}^{p}(G:F).$ For $p=2$ and $G=SL(2,\mathbb{R}),$ no such function exists and we have (in this example) that $$\bigoplus_{F\subset\hat{K},\\|F|<\infty}\mathcal{C}^{2}(G:F)=\bigoplus_{m,n\in\mathbb{Z}}\mathcal{C}^{2}_{mn}(G)=\mathcal{C}^{2}_{mn}(G),$$ Varadarajan $[24.],$ p. $280-281.$ Indeed, a readily available finite set $F$ for a clear understanding of the proof of Proposition $4.2$ is $F=\{\xi_{1},\xi_{2}\}$ for any $\xi_{1},\xi_{1}\in\hat{K}.$
	
	\indent $(2)$ The proof of Proposition $4.2$ shows that every function $\alpha\in\mathcal{C}^{p}(G)$ has a unique expansion into the sum $$\alpha=\sum_{F\subset\hat{K},\\|F|<\infty}\bar{\chi}_{_{F}}*\alpha*\bar{\chi}_{_{F}}$$ in $\bigoplus_{F\subset\hat{K},\\|F|<\infty}\mathcal{C}^{p}(G:F),$ where we may denote the convolution $\bar{\chi}_{_{F}}*\alpha*\bar{\chi}_{_{F}}$ simply by $\alpha_{F}.$ Now as each $\mathcal{C}^{p}(G:F)$ splits into a direct sum given explicitly as $$\mathcal{C}^{p}_{B}(G:F)\oplus\mathcal{C}^{p}_{H}(G:F)$$ (if $rk(G)=rk(K)$) or as $\mathcal{C}^{p}_{H}(G:F)$ (if $rk(G)>rk(K)$), it follows that $\alpha_{F}$ is a sum which may be given as $\alpha_{F}=\alpha_{B,F}+\alpha_{H,F},$ where $\alpha_{B,F}\in\mathcal{C}^{p}_{B}(G:F)$ and $\alpha_{H,F}\in\mathcal{C}^{p}_{H}(G:F).$ Hence, each $\alpha\in\mathcal{C}^{p}(G)$ may be uniquely written as $\alpha=\alpha_{B}+\alpha_{H},$ with $$\alpha_{B}=\sum_{F\subset\hat{K},\\|F|<\infty}\bar{\chi}_{_{F}}*\alpha_{B}*\bar{\chi}_{_{F}}$$ and $$\alpha_{H}=\sum_{F\subset\hat{K},\\|F|<\infty}\bar{\chi}_{_{F}}*\alpha_{H}*\bar{\chi}_{_{F}}.$$\\
	
	\indent {\bf 4.4 Corollary.} Let $F\subset\hat{K}$ with $|F|<\infty.$ Then $$\bigoplus_{F\subset\hat{K},|F|<\infty}\mathcal{C}^{p}(G:F)$$ is a dense Fr$\acute{e}$chet convolution Schwartz subalgebra of $\mathcal{C}^{p}(G).\;\Box$\\
	
	\indent We may now consider the (tensor) product of Trombi's images of the Fourier transform of each of $\mathcal{C}^{p}(G:F),$ as suggested by the $SL(2,\mathbb{R})-$versions of  Barker $[5.]$ and Varadarajan $[24].$\\
	
	\indent {\bf 4.5 Proposition.} Let $F\subset\hat{K}$ with $|F|<\infty.$ Then $$\bigotimes_{F\subset\hat{K},|F|<\infty}\mathcal{C}^{p}(\hat{G}:F)$$ is a dense Fr$\acute{e}$chet multiplication Schwartz subalgebra of $\mathcal{C}^{p}(\hat{G}).$
	
	\indent {\bf Proof.} Each of $\mathfrak{F}^{F}:\mathcal{C}^{p}(G:F)\rightarrow\mathcal{C}^{p}(\hat{G}:F)$ is a linear topological algebra isomorphism, Trombi $[20.],$ p. $122.$ We now define the map $\mathfrak{F}$ on $\bigoplus_{F\subset\hat{K},|F|<\infty}\mathcal{C}^{p}(G:F)$ as $$\mathfrak{F}(\alpha_{1},\alpha_{2},\cdots,\alpha_{n},\cdots):=
	\mathfrak{F}^{F_{1}}(\breve{\alpha}_{1})\otimes\mathfrak{F}^{F_{2}}(\breve{\alpha}_{2})\otimes\cdots\otimes\mathfrak{F}^{F_{n}}(\breve{\alpha}_{n})\otimes\cdots,$$ where $F_{j}\subset\hat{K},$ $|F_{j}|<\infty$ and $\alpha_{j}\in\mathcal{C}^{p}(G:F_{j}).$ See Barker $[5.],$ p. $44$ for the $SL(2,\mathbb{R})$ version of $\mathfrak{F}.$ The result follows from Corollary $4.4,$ since each of $$\mathfrak{F}^{F}:\mathcal{C}^{p}(G:F)\rightarrow\mathcal{C}^{p}(\hat{G}:F)$$ is a linear topological algebra isomorphism$.\;\Box$\\
	
	\indent We are now in a position to state the fundamental theorem of harmonic analysis on any $G$ of real-rank $1$ via the map $\mathfrak{F}$ in the proof of Proposition $4.5.$\\
	
	\indent {\bf 4.6 Theorem.} The Fourier transform $$\mathfrak{F}:\mathcal{C}^{p}(G)\rightarrow\mathcal{C}^{p}(\hat{G}):\alpha\mapsto\mathfrak{F}(\breve{\alpha})$$ is a linear topological algebra isomorphism of the Schwartz Fr$\acute{e}$chet algebras.
	
	\indent {\bf Proof.} The result is an immediate consequence of Corollary $4.4$ and Proposition $4.5.\;\Box$\\
	
	\indent {\bf 4.7 Remarks.}\\
	\indent $(1)$ The canonical wave-packet corresponding to every $\alpha\in\mathcal{C}^{p}(G)$ in the proof of Theorem $4.6$ is here denoted by $\phi_{_{\mathfrak{F}(\alpha)}}$ and is the function on $G$ given as $$\phi_{_{\mathfrak{F}(\alpha)}}(x)=\sum_{\Lambda\in L^{+}_{B}(F)}d(\Lambda)\phi_{_{\mathfrak{F}_{B}(\alpha)}}(\Lambda:x)+\sum_{\chi\in \hat{M}(F)}d(\chi)\phi_{_{\mathfrak{F}_{H}(\alpha)}}(Q:\chi:x),$$ where $$\phi_{_{\mathfrak{F}_{B}(\alpha)}}(\Lambda:x)=
	\sum_{F\subset\hat{K},|F|<\infty}tr\{\mathfrak{F}^{F}_{B}(\alpha)(\Lambda)\pi^{F}_{\Lambda}(x)\}$$ and $$\phi_{_{\mathfrak{F}_{H}(\alpha)}}(Q:\chi:x)=
	\sum_{F\subset\hat{K},|F|<\infty}tr\{\mathfrak{F}^{F}_{H}(\alpha)(Q:\chi:\nu)\pi^{F}_{Q:\chi:\nu}(x)\}\mu^{F}(\chi:\nu)$$ are the corresponding wave-packets in $\mathcal{C}^{p}_{B}(G:F)$ and $\mathcal{C}^{p}_{H}(G:F),$ respectively. (See Trombi $[20.],$ p. $85$). See also Barker $[5.],$ p. $49-50.$\\
	\indent $(2)$ The construction of the dense Fr$\acute{e}$chet convolution Schwartz subalgebra in Corollary $4.4$ is in line with the splitting of $\mathcal{C}^{p}(G)$ in Barker $[5.],$ p. $44$ for $G=SL(2,\mathbb{R}).$ It thus follows that the Fourier transform $\mathfrak{F}$ maps a direct sum in $\mathcal{C}^{p}(G)$ (on)to a tensor product in $\mathcal{C}^{p}(\hat{G}).$ This suggests a construction of a basis in both $\mathcal{C}^{p}(G)$ and $\mathcal{C}^{p}(\hat{G}).$
	
	\indent A consideration of the basis $$\phi_{_{\gamma,m}}\otimes\phi_{_{\gamma,l}},$$ with $\gamma\in F$ and $(l,m)\in\mathbb{Z}^{2}$ of the space $$\mathcal{C}^{p}(\hat{G}:F) = \begin{cases}
		\;\;\;\;\;\;\;\;\mathcal{C}^{p}_{H}(\hat{G}:F),  \;\;\;\;\;\;\;\;\;\;\;\;\;\;\;\;\;\text{if}\;\;\; rk(G)>rk(K), \\
		\mathcal{C}^{p}_{B}(\hat{G}:F)\otimes\mathcal{C}^{p}_{H}(\hat{G}:F),  \;\;\;\;\;\text{if}\;\;\; rk(G)=rk(K),
	\end{cases}$$ reveals that the members of $\mathcal{C}^{p}(\hat{G}:F)$ could be realized as matrices of the form $$(\mathfrak{F}_{B}(\breve{\alpha})_{(\gamma,m)}(\Lambda)\otimes\mathfrak{F}_{H}(\breve{\alpha})_{(\gamma,l)}(Q:\chi:\nu))_{\gamma,M},$$ where $\gamma\in F$ and $M:=(l,m)\in\mathbb{Z}^{2}.$ We can therefore explicitly realize $\mathcal{C}^{p}(\hat{G})$ as an algebra of infinite matrices with a countable number of blocks as follows.\\
	
	\indent {\bf 4.8 Corollary.} The Fourier image $\mathcal{C}^{p}(\hat{G})$ (of $\mathcal{C}^{p}(G)$ under the map $\mathfrak{F}$) is a Schwartz Fr$\acute{e}$chet multiplication algebra consisting of block matrices of the form $$((\mathfrak{F}_{B}(\breve{\alpha})_{(\gamma,m)}(\Lambda)\otimes\mathfrak{F}_{H}(\breve{\alpha})_{(\gamma,l)}(Q:\chi:\nu))_{\gamma\in F, (l,m)\in\mathbb{Z}^{2}})_{F\subset \hat{K},|F|<\infty}.$$
	
	\indent {\bf Proof.} The matrix realization follows from combining Corollary $3.3,$ Proposition $4.5,$ the above remark on a basis of $\mathcal{C}^{p}(\hat{G}:F)$ and Theorem $4.6.\;\Box$\\
	
	\indent The reader may consult Varadarajan $[24.],$ p. $276-278$ for the Hilbert-space version of the last corollary when $p=2$ and $G=SL(2,\mathbb{R}).$ In this special (linear) case, the block-matrix of Corollary $4.8$ reduces (or is isomorphic) to the partitioned block-matrix $$((\int f\bar{\phi}_{kmn}dG),(\int f\bar{f}_{mn}(\mu:\cdot)dG))_{m,n\in\mathbb{Z}}$$ for $f\in\mathcal{C}^{2}(SL(2,\mathbb{R}))$ as contained in Varadarajan $[24.].$
	
	\indent Now if $\mathcal{C}^{p}(G)^{'}$ and $\mathcal{C}^{p}(\hat{G})^{'}$ denote the topological dual spaces of $\mathcal{C}^{p}(G)$ and $\mathcal{C}^{p}(\hat{G}),$ respectively, with the weak topology, the transpose $\mathfrak{F}^{'}$ (See the end-matters of Arthur $[1]$) of $\mathfrak{F}$ in Theorem $4.4$ exists and is equally a linear topological isomorphism $\mathcal{C}^{p}(\hat{G})^{'}\rightarrow\mathcal{C}^{p}(G)^{'}.$ This isomorphism affords the development of harmonic analysis for $p-$tempered distributions on groups $G$ of real-rank $1.$ See also Oyadare $[16.]$ and $[17].$
	
	\indent The success of the methods used above hinges on the earlier established fact that almost all members of the space $\mathcal{C}^{p}(G)$ may be realized as a countable sum of some $\tau=(\tau_{l},\tau_{r})-$spherical functions on $G.$ This allowed us to use some density properties of the space of these functions to arrive at our conclusion for real-rank $1$ groups $G.$ These same density properties are actually available for any $G$ of any arbitrary rank and we shall explore the existence of such density results in arriving at a theorem similar to Theorem $4.6$ for such $G$ via the space of $K-$finite functions in $\mathcal{C}^{p}(G)$ in the next section.\\
	\ \\
	\S5. {\bf Arbitrary real-rank case for $G.$}
	
	\indent Now let $G$ denote any real reductive group of arbitrary real-rank. We shall refer to Varadarajan $[25.]$ for the structure theory of $G.$ The method employed for the arbitrary real-rank case is to combine the techniques of Varadarajan $[24.]$ and G$\dot{a}$rding $[12.].$ We recall that for any representation $\pi$ of $G$ on a vector space $V,$ the first known proof (if not the only one still known) that verifies that the subspace $C^{\infty}(\pi)$ (of $V$) consisting of smooth vectors for $\pi$ is dense in $V,$ is indirectly through the use of what is now known as the G$\dot{a}$rding subspace $\mathcal{G}(\pi).$ The indirect proof (of the denseness of $C^{\infty}(\pi)$ in $V$) was achieved in two simple steps:\\
	$$\mathcal{G}(\pi)\subseteq C^{\infty}(\pi)\subseteq V\;\;\;\;\mbox{and}\;\;\;\; \mathcal{G}(\pi)\;\mbox{is dense in}\;V.$$
	
	\indent The argument being put forward in this paper is analogous to the method of L. G$\dot{a}$rding $[12.].$ Indeed, in its technical details, it is the same as using the approach of type$-(m,n)$ spherical functions as done by Varadarajan $[24.]$ to create subspaces whose direct sum in $\mathcal{C}^{p}(G)$ (respectively, tensor product in $\mathcal{C}^{p}(\hat{G})$) is dense in $\mathcal{C}^{p}(G)$ (respectively, in $\mathcal{C}^{p}(\hat{G})$). The first of these series of results is the following.\\
	
	\indent {\bf 5.1 Proposition.} We have that $$\mathcal{C}^{p}_{\tau_{l},\tau_{r}}(G)\subseteq\mathcal{C}^{p}(G:F)\subseteq\mathcal{C}^{p}(G),$$ for any $\tau_{l},\tau_{r}$ in some $F$ in $\hat{K}.$
	
	\indent {\bf Proof.} It is clear that $\{\tau_{l},\tau_{r}\}$ is a finite subset of $\hat{K}$ on which the Trombi analysis of $[20.]$ restricts to$.\;\Box$\\
	
	\indent {\bf 5.2 Theorem.} The direct sum $\bigoplus_{\tau_{l},\tau_{r}\in\hat{K}}\mathcal{C}^{p}_{\tau_{l},\tau_{r}}(G)$ is dense in $\mathcal{C}^{p}(G).$
	
	\indent {\bf Proof.} Every $f\in\mathcal{C}^{p}(G)$ has a representation as $$f=\sum_{F\subset\hat{K},\\|F|<\infty}\bar{\chi}_{_{F}}*f*\bar{\chi}_{_{F}},$$ where each pair $\{\tau_{l},\tau_{r}\}\subset F.\;\Box$\\
	
	\indent {\bf 5.3 Theorem.} The direct sum $\bigoplus_{F\subseteq\hat{K},|F|<\infty}\mathcal{C}^{p}(G:F)$ is dense in $\mathcal{C}^{p}(G).$
	
	\indent {\bf Proof.} This follows from Theorem $5.2.\;\Box$\\
	
	\indent Now denote by $\mathcal{C}^{p}_{\tau_{l},\tau_{r}}(\hat{G}),$ $\mathcal{C}^{p}(\hat{G}:F)$ and $\mathcal{C}^{p}(\hat{G})$ the respective images of $\mathcal{C}^{p}_{\tau_{l},\tau_{r}}(G),$ $\mathcal{C}^{p}(G:F)$ and $\mathcal{C}^{p}(G)$ under the operator-valued Fourier transform $\alpha\mapsto\mathfrak{F}(\breve{\alpha}).$ The fact that $\mathfrak{F}$ is a linear topological isomorphism given as $$\mathcal{C}^{p}_{\tau_{l},\tau_{r}}(G)\rightarrow\mathcal{C}^{p}_{\tau_{l},\tau_{r}}(\hat{G})$$ and $$\mathcal{C}^{p}(G:F)\rightarrow\mathcal{C}^{p}(\hat{G}:F)$$ in Trombi $[20.]$ imply the following result.\\
	
	\indent {\bf 5.4 Theorem.} We always have that $$\mathcal{C}^{p}_{\tau_{l},\tau_{r}}(\hat{G})\subseteq\mathcal{C}^{p}(\hat{G}:F)\subseteq\mathcal{C}^{p}(\hat{G}),$$ for any $\tau_{l},\tau_{r}$ in some $F$ and that $\bigotimes_{\tau_{l},\tau_{r}\in\hat{K}}\mathcal{C}^{p}_{\tau_{l},\tau_{r}}(\hat{G})$ is dense in $\mathcal{C}^{p}(\hat{G}).$ Hence $\bigotimes_{F\subset\hat{K},|F|<\infty}\mathcal{C}^{p}(\hat{G}:F)$ is dense in $\mathcal{C}^{p}(\hat{G}).\;\Box$\\
	
	\indent We then have the full version of Theorem $4.6$ given below. \\
	
	\indent {\bf 5.5 Theorem.} The Fourier transform $$\mathfrak{F}:\mathcal{C}^{p}(G)\rightarrow\mathcal{C}^{p}(\hat{G}):\alpha\mapsto\mathfrak{F}(\breve{\alpha})$$ is a linear topological algebra isomorphism of the Schwartz Fr$\acute{e}$chet algebras$.\;\Box$\\
	
	\indent This leads to an explicit expression for members of $\mathcal{C}^{p}(\hat{G})$ in terms of the earlier basis given in \S 4. This is explicitly given below.\\
	
	\indent {\bf 5.6 Corollary.} The Fourier image $\mathcal{C}^{p}(\hat{G})$ (of $\mathcal{C}^{p}(G)$ under the map $\mathfrak{F}$) is a Schwartz Fr$\acute{e}$chet multiplication algebra consisting of block matrices of the form $$((\mathfrak{F}_{B}(\breve{\alpha})_{(\gamma,m)}(\Lambda)\otimes\mathfrak{F}_{H}(\breve{\alpha})_{(\gamma,l)}(Q:\chi:\nu))_{\gamma\in F, (l,m)\in\mathbb{Z}^{2}})_{F\subset \hat{K},|F|<\infty}.\;\Box$$\\
	
	\S6. {\bf Coda}
	
	\indent Another approach to the general theory is the following, based on Trombi's $K-$finite results. Consider the Fr$\acute{e}$chet convolution Schwartz algebra $\mathcal{C}^{p}(G,K)$ consisting of $K-$finite members of $\mathcal{C}^{p}(G)$. It has been shown by Trombi (See Barker $[5.]$) that $\mathfrak{F}:\mathcal{C}^{p}(G,K)\rightarrow\mathcal{C}^{p}(\hat{G},K)$ is a linear topological algebra isomorphism of Fr$\acute{e}$chet algebras. Let us now define a subalgebra of $\mathcal{C}^{p}(G)$ denoted as $\mathcal{C}^{p}(G,K:F)$ explicitly as $$\mathcal{C}^{p}(G,K:F):=\mathcal{C}^{p}(G,K)\bigcap\mathcal{C}^{p}(G:F).$$ This is the subalgebra of all $K-$finite members of $\mathcal{C}^{p}(G)$ whose $K-$types lie in a finite subset $F$ of $\hat{K}.$ $\mathcal{C}^{p}(G,K:F)$ is topologized by the induced continuous seminorms of $\mathcal{C}^{p}(G,K).$ We equally define $\mathcal{C}^{p}_{H}(G,K:F),$ $\mathcal{C}^{p}_{H}(G,K:F)$ and $\mathcal{C}^{p}(G,K:F)$ as in the earlier case of real-rank $1.$ We have the decomposition of $\mathcal{C}^{p}(G,K:F)$ (respectively, $\mathcal{C}^{p}(\hat{G},K:F):=\mathfrak{F}(\mathcal{C}^{p}(G,K:F))$) into the direct sum (respectively, tensor product) of the space $\mathcal{C}^{p}_{H}(G,K:F)$ with $\mathcal{C}^{p}_{B}(G,K:F)$ or with $0$ (respectively, $\mathcal{C}^{p}_{H}(\hat{G},K:F)$ with $\mathcal{C}^{p}_{B}(\hat{G},K:F)$ or with $0$), depending on whether $rk(G)=rk(K)$ or $rk(G)>rk(K),$ respectively. It is clear that $$\bigoplus_{F\subset\hat{K},|F|<\infty}\mathcal{C}^{p}(G,K:F)=\mathcal{C}^{p}(G,K),$$ which by Trombi (in Barker $[5.]$) is dense in $\mathcal{C}^{p}(G).$ We then have a more general version of Proposition $4.5$ in that for $F\subset\hat{K}$ with $|F|<\infty,$ $$\bigotimes_{F\subset\hat{K},|F|<\infty}\mathcal{C}^{p}(\hat{G},K:F)$$ is the Fr$\acute{e}$chet multiplication algebra $\mathcal{C}^{p}(\hat{G},K)$ and is dense in $\mathcal{C}^{p}(\hat{G}).$ This also leads to the fundamental theorem on $G$ as follows.\\
	
	\indent {\bf 6.1 Theorem.}  The Fourier transform $$\mathfrak{F}:\mathcal{C}^{p}(G)\rightarrow\mathcal{C}^{p}(\hat{G}):\alpha\mapsto\mathfrak{F}(\breve{\alpha})$$ is a linear topological algebra isomorphism of the Fr$\acute{e}$chet convolution Schwartz algebras$.$
	
	\indent {\bf Proof.} Each $\mathfrak{F}:\mathcal{C}^{p}(G,K:F)\rightarrow\mathcal{C}^{p}(\hat{G},K:F)$ is a linear topological algebra  isomorphism. Since $\bigoplus_{F\subset\hat{K},|F|<\infty}\mathcal{C}^{p}(G,K:F)=\mathcal{C}^{p}(G,K)$ is dense in $\mathcal{C}^{p}(G)$ and $\bigotimes_{F\subset\hat{K},|F|<\infty}\mathcal{C}^{p}(\hat{G},K:F)=\mathcal{C}^{p}(\hat{G},K)$ is dense in $\mathcal{C}^{p}(\hat{G}),$ then the result. follows$.\;\Box$\\
	
	\indent The proof of the Bochner theorem and analysis of its canonical wave-packets on $G,$ corresponding to the generalization contained in both Theorems $4.6$ (for real-rank $1$ groups $G$) and Theorem ($5.5$ or) $6.1$ (for arbitrary real-rank groups $G$) is essentially as established in Oyadare $[17.]$ $[18.]$ with necessary basis-adjustment (from the scalar version to the vector version) as required.\\

	\ \\
	{\bf References.}
	\begin{description}
		\item [{[1.]}] Arthur, J. G., \textit{Harmonic analysis of tempered distributions on semisimple Lie groups of real rank one,} Ph.D. Dissertation, Yale University, $1970.$
		
		\item [{[2.]}] Arthur, J. G., \textit{Harmonic analysis of the Schwartz space of a reductive Lie group I,} mimeographed note, Yale University, Mathematics Department, New Haven, Conn;
		
		\item [{[3.]}] Arthur, J. G., \textit{Harmonic analysis of the Schwartz space of a reductive Lie group II,} mimeographed note, Yale University, Mathematics Department, New Haven, Conn.
		
		\item [{[4.]}] Bargmann, V.,  Irreducible unitary representations of the Lorentz group, {\it Ann. of Math.} vol. {\bf 48}, $(1947),$ p. $568-640..$
		
		\item [{[5.]}] Barker, W. H.,  $L^{p}$ harmonic analysis of $SL(2,\mathbb{R}),$ \textit{American Mathematical Society Memoirs,} vol. \textbf{76} No. \textbf{393}, ($1988$).
		
		\item [{[6.]}] Dixmier, J., Op$\acute{\mbox{e}}$rateurs de rang fini dans les repr$\acute{\mbox{e}}$sentations unitaires,\textit{ Publ. math. de l' Inst. Hautes $\acute{\mbox{E}}$tudes Scient.,} tome $\textbf{6}$ ($1960$), p. $13-25.$
		
		\item [{[7.]}] Eguchi, M., The Fourier transform of the Schwartz space on a semisimple Lie group. {\it Hiroshima Math. J.}, vol. {\bf 4}, $(1974),$ p. $133-209.$
		
		\item [{[8.]}] Ehrenpreis L. and Mautner F., Some properties of the Fourier transform on semisimple Lie groups, I. {\it Ann. Math.,} {\bf 61,} $(1955),$ p. $406-439.$
		
		\item [{[9.]}] Ehrenpreis L. and Mautner F., Some properties of the Fourier transform on semisimple Lie groups, II. {\it Trans. Amer. Math. Soc.}, vol. {\bf 84}, $(1957),$ p. $1-55.$
		
		\item [{[10.]}] Ehrenpreis L. and Mautner F., Some properties of the Fourier transform on semisimple Lie groups, III. {\it Trans. Amer. Math. Soc.}, vol. {\bf 90}, $(1959),$ p. $431-484.$
		
		\item [{[11.]}] Gangolli, R. and Varadarajan, V. S., \textit{Harmonic analysis of spherical functions on real reductive groups,} Ergebnisse der Mathematik und iher Genzgebiete, $vol.$ {\bf 101}, Springer-Verlag, Berlin-Heidelberg. $1988.$
		
		\item [{[12.]}] G$\dot{a}$rding L., Note on continuous representations of Lie groups, {\it Proc. Nat. Acad. Sci. USA,} vol. {\bf 33,} $(1947),$ p. $331-332.$
		
		\item [{[13.]}] Harish-Chandra, Harmonic analysis on semisimple Lie groups, {\it Bull. AMS,} vol. {\bf 78,} $(1970),$ p. $529-551.$
		
		\item [{[14.]}] Knapp, A. W., {\it Representation theory os semisimple Groups: An Overview Based on Examples,} Princeton University Press, Princeton, $(1986).$ Reprinted: $2001.$
		
		\item [{[15.]}] Langlands, R. P., Unpublished manuscript. Reported by Harish-Chandra in $[13.]$
		
		\item [{[16.]}] Oyadare, O. O.,  Non-spherical Harish-Chandra Fourier transforms on real reductive groups, \textit{J. Fourier Anal. Appl.} \textbf{28}, $15$ $(2022).$\\ http://doi.org/10.1007/s00041-09906-w
		
		\item [{[17.]}] Oyadare, O. O.,  The full Bochner theorem on real reductive groups, \textit{Algebras, Groups and Geometries} \textbf{39}, $(2023),$ p. $207-220.$
		
		\item [{[18.]}] Oyadare, O. O.,  Functional analysis of canonical wave-packets on real reductive groups, \textit{arXiv:$1912.07542v1.$} [math.FA], $13$ Dec. $2019.$
		
		\item [{[19.]}] Trombi, P. C., Spherical transforms on symmetric spaces of rank	one (or Fourier analysis on semisimple Lie groups of split rank one), {\it Thesis, University of Illinios} $(1970).$
		
		\item [{[20.]}] Trombi, P. C., Harmonic analysis of $\mathcal{C}^{p}(G : F),\; (1 \leq p < 2)$ {\it J. Funct. Anal.,} vol.	{\bf 40.} $(1981),$ p. $84-125.$
		
		\item [{[21.]}] Trombi, P. C., Invariant harmonic analysis on split rank one groups with applications. {\it Pacific J. Math.} vol. {\bf 101.} no.: {\bf 1.} $(1982),$ p. $223-245.$
		
		\item [{[22.]}] Trombi, P. C., Unpublished manuscript. Reported by W. H. Barker in $[5.]$
		
		\item [{[23.]}] Trombi, P. C. and Varadarajan, V. S., Spherical transforms on semisimple
		Lie groups, {\it Ann. Math.,} vol. {\bf 94.} $(1971),$ p. $246-303.$
		
		\item [{[24.]}] Varadarajan, V. S., {\it An introduction to harmonic analysis on semisimple
			Lie groups,} Cambridge Studies in Advanced Mathematics, {\bf 161,} Cambridge University Press, $1989.$
		
		\item [{[25.]}] Varadarajan, V. S., {\it Harmonic analysis on real reductive groups,} Lecture Notes in Mathematics, {\bf 576,} Springer Verlag, $1977.$
		
		\item [{[26.]}] Wallach, N.,  {\it Harmonic analysis in homogeneous spaces,} Dekker, New York, $1973.$
		
		\item [{[27.]}] Warner, G., {\it Harmonic analysis on semisimple Lie groups, I.} Springer-Verlag, New York, $1972.$
		
	\end{description}
	
\end{document}